\documentclass{birkjour}
 \usepackage{amsmath, amsfonts, amssymb, amsthm, euscript, amscd, latexsym, mathrsfs}
 \usepackage{lineno,hyperref}

\newtheorem{thm}{\sc Theorem}
\newtheorem{lem}{\sc Lemma}
\newtheorem{cor}{\sc Corollary}

\newtheorem{rem}{\sc Remark}

\newcommand{\rf}{\eqref}
\newcommand{\sss}{\scriptscriptstyle}

\newcommand{\eps}{\varepsilon}
\newcommand{\pl}{\partial}
\newcommand{\gt}{\geqslant}
\newcommand{\lt}{\leqslant}
\newcommand{\sub}{\subset}

\newcommand{\dl}{\delta}
\newcommand{\lb}{\label}
\newcommand{\al}{\alpha}
\newcommand{\gm}{\gamma}
 
 \newcommand{\Dl}{\Delta}
 \newcommand{\la}{\lambda}
 \newcommand{\sg}{\sigma}
\newcommand{\dd}{\diagdown}
\newcommand{\diag}{{\rm \, diag \,}}

\newcommand{\mc}{\mathcal}

\newcommand{\mb}{\mathbf}

\newcommand{\td}{\tilde}

\newcommand{\x}{\times}
\newcommand{\mto}{\mapsto}

\newcommand{\C}{{\rm C}}

\newcommand{\Te}{\Theta}
\newcommand{\te}{\theta}

\newcommand{\lap}{\Delta}
\newcommand{\nab}{\nabla}

\newcommand{\fdot}{\,\cdot\,}
\newcommand{\scal}{{\rm scal}}
\def\Rnu{{\mathbb R}}

\def\com#1{}

\def\aa#1{ \begin{eqnarray*} #1 \end{eqnarray*} }
\def\aaa#1{ \begin{eqnarray} #1 \end{eqnarray} }
\def\mm#1{ \begin{multline*} #1 \end{multline*} }
\def\mmm#1{ \begin{multline} #1 \end{multline} }
\newtheorem*{thm*}{\sc Theorem}

\begin{document}

\title[Uniform approximation of the heat kernel on a manifold]{Uniform approximation of the heat kernel on a manifold}
\author{Evelina Shamarova} 
\address{Mathematics department, Federal University of Para\'iba, 58051, Jo\~ao Pessoa, Brazil}
\email{evelina@mat.ufpb.br}
\author{Alexandre B. Simas}
\address{Mathematics department, Federal University of Para\'iba, 58051, Jo\~ao Pessoa, Brazil}
\email{alexandre@mat.ufpb.br}

 \begin{abstract}
We approximate the heat kernel $h(x,y,t)$ on a compact connected Riemannian manifold $M$ without boundary uniformly in $(x,y,t)\in M\times M\times [a,b]$, $a>0$, by $n$-fold integrals over $M^n$ of the densities of Brownian bridges.
Moreover, we provide an estimate for the uniform convergence rate.
 As an immediate corollary, we get a uniform approximation of solutions of the Cauchy problem for the heat equation on $M$.
   \end{abstract}

\thanks{
The research of A.B. Simas was funded by 
the National Counsel of Technological and Scientific Development (CNPq, Brazil).}

\keywords{Heat kernel, Heat equation, Riemannian manifold} 
\subjclass{35K08, 58J35}

\maketitle

 \section{Introduction}

 The heat kernel  on a manifold is widely used in diverse areas of mathematics and mathematical physics,
 in particular, geometric analysis and quantum field theory.
 It is defined as the fundamental solution of the \textit{heat equation} on the manifold $M$
 \aaa{
 \lb{heat}
  \pl_t u =-\frac12 \lap_M u, 
  }
 where 
 $\lap_M$ is the Laplace-Beltrami operator. 
 In the past decades the so called heat kernel approach 
as well as various heat kernel estimates and asymptotic expansions arose  in
 mathematics and physics literature. 
 We refer the reader to \cite{grigoryan}  for the extensive bibliography on the heat kernel
 estimates, methods, and applications.

In this work we represent the heat kernel $h(x,y,t)$ on $M$ as a uniform limit
 of transition density functions of Brownian bridge type processes in a Euclidean space $\Rnu^m$. 
Below we list possible applications of our formula. It allows one to uniformly approximate the solution of the Cauchy problem
for \rf{heat} while maintaining control of the error.
Furthermore, our formula allows  to construct computational schemes for the heat kernel.
 It involves integration over $M$ which, under some conditions, can be approximately performed even when $M$ 
is not known, but instead, we have an approximation of $M$ by a mesh. Finally, the uniform approximation
of the heat kernel is important because of the absence of an explicit formula for $h(x,y,t)$ for most of manifolds.

Let $M$ be a $d$-dimensional compact connected Riemannian manifold without boundary
 isometrically embedded into $\Rnu^m$, $m\gt d+1$, by the Nash theorem (\cite{nash}).
Define 
\aaa{
\lb{pdens}
p(x,y,t)= (2\pi t)^{-\frac{d}2} e^{-\frac{|x-y|^2}{2t}}
\; 
\text{and} \; \;  q(x,y,t)= \frac{p(x,y,t)}{\int_M p(x,z,t) \la_M(dz)},
} 
 where  $\la_M$ is the volume measure on $M$ and $|\cdot|$ is the Euclidean norm in $\Rnu^m$.
 Let  $\mc P = \{ 0=t_0 < t_1 < \cdots < t_{n+1}=t \}$ be a partition
 of the interval $[0,t]$, and let $|\mc P|$ denote the mesh of $\mc P$. Define
 \mmm{
 \lb{qp0}
 q_{\sss \mc P}(x, y,t) = 
  \int_M \hspace{-1mm} \la_M(dx_1) \, q(x, x_1, \Dl t_0) \\
  \cdots 
   \int_M  \hspace{-1mm} \la_M(dx_n)  q(x_{n-1}, x_{n}, \Dl t_{n-1})\,
  q(x_{n}, y, \Dl t_{n}),
  }
  where $\Dl t_i = t_{i+1} - t_{i}$, $i=0, \ldots, n$.
  Below we state our main result.
\begin{thm}[Uniform approximation of the heat kernel]
\lb{main}
Let $M$ be a compact connected Riemannian manifold without boundary.
Assume there exist a number $L>1$ and an integer $n_0>1$ so that  for all $n\gt n_0$  it holds that
 $\min_{0\lt i \lt n} \Dl t_i > |\mc P|^L$. Then, as  $|\mc P|\to 0$,
 \aaa{
 \lb{conv-main}
q_{\sss \mc P}(x, y, t) \to h(x, y, t)
}
uniformly in $(x,y,t)\in M \x M \x [a,b]$, where $a>0$.
 \end{thm}
 Additionally, in Theorem \ref{rate} we prove that the rate of  uniform
  convergence \rf{conv-main} is of order $\mathcal{O}(|\mc P|^{\frac1{2(d+9)}})$.
 
 Our article was inspired by the work of Smolyanov et al.  \cite{sm_vr}
about the weak approximation of the Wiener measure on paths on $M$ 
by the distributions of
successive Brownian bridges in $\Rnu^m$ returning to $M$ at each time $t_i$
whose transition density functions are $q_{\sss \mc P}(x, y,t)$.

 For the proof of Theorem \ref{main}, we establish relations between $h(x,y,t)$ and $q(x,y,t)$ 
 in a $t^\al$-neighborhood of the diagonal of $M\x M$,
 where $t$ is a small parameter.
Also, we make use of 
 an analytical result of  Smolyanov et al. \cite{sm_vr} (see Corollary 2 on p. 596).
 We apply this result to the density $p(x,y,\tau)$ which depends on another small parameter $\tau$.
 Due to the dependence on $\tau$,  the application of the result of Smolyanov et al. is not straightforward
 but requires a delicate analysis to choose the parameters $\tau$ and $\al$
 to guarantee the uniform convergence.
 We would like to remark that although it was not the goal of the article of Smolyanov et al. \cite{sm_vr} 
 to investigate convergence \rf{conv-main}, the results of this work imply
 \rf{conv-main} weakly in $y\in M$ and pointwise in $(x,t)\in M\x (0,\infty)$, which, in turn,
 implies the pointwise approximation of solutions of the Cauchy problem for \rf{heat}. 
 Hence, the result of this paper can be regarded as an improvement with respect to the aforementioned 
 approximation since ours is uniform for both, the heat kernel
 and the solution of the Cauchy problem for  \rf{heat}.
 \section{Theorems on uniform approximation}
 In this section, we will prove Theorem \ref{main} on uniform approximation of the heat kernel on a manifold and will give an estimate
 of the convergence rate via the mesh of the partition $\mc P$. Moreover, we derive an approximation formula for the heat kernel on the $d$-dimensional sphere
 $S^d$ in $\Rnu^{d+1}$ when the partition $\mc P$ is uniform.
 
 For $x,y\in M$, $t>0$,
 we define
 $\mc E(x,y,t) = (2\pi t)^{-\frac{d}2} e^{-\frac{\rho(x,y)^2}{2t}}$,
 where $\rho(x,y)$ is the geodesic distance. By the Nash theorem,
  $M$ is assumed to be isometrically embedded into $\Rnu^m$. 
 Let for any $\eps>0$, $U_\eps(\diag (M\x M))$ denote the $\eps$-neighborhood
  of the diagonal of $M\x M$, i.e. $ \{(x,y) \in M\x M: |x-y| < \eps\}$.
  Define the set
  \aa{
  D_{\eps,\al} = \{(x,y,t): t\in (0,\eps), (x,y) \in U_{t^\al}(\diag (M\x M))\}.
  }

 \subsection{Proof of Theorem \ref{main}}
 For the proof of the theorem we will need a few lemmas.
 \begin{lem}
 \label{lem1}
Let $\al \in (\frac14, \frac12)$. Then, there exists $\eps>0$, 
and functions $\Te, R: D_{\eps,\al} \to \Rnu$ bounded uniformly in $\al$ and
$(t,x,y)\in D_{\eps,\al}$, so that on $D_{\eps,\al}$ 
\aaa{
 \label{2nd}
  q(x,y,t)  =  h(x,y,t)(1+\Te(x,y,t) \,t^{4\al-1}) + R(x,y,t) t^3.
   }
  \end{lem}
 \begin{proof}
For a function $f\in \C^2(M)$, recall the asymptotic expansion (see \cite{shamarova})
\mmm{
\lb{expansion}
\int_M f(y) p(x,y,t) \la_M(dy)
   =  f(x) - f(x) t \Big( \frac16 \scal(x)
   +   \frac1{16} \lap_M\lap_M\left.|x-\,\cdot\,|^2\right|_{x}\Big) \\
   -\frac{t}2 \lap_M f(x)
   +  t^{\frac32} \mc R(x,t),
}
 where $\mc R(x,t)$ is bounded in both arguments,
 and $\scal(x)$ is the scalar curvature at $x$.
Applying \rf{expansion} to $f=1$, we obtain
that there exist $\eps_1>0$ and a bounded function $\Te: M\x M\x (0,\eps_1)\to \Rnu$, so that
that for all $x,y\in M$ and $t<\eps_1$,
\aaa{
\lb{pq}
q(x,y,t)  =  p(x,y,t)(1+ \Te(x,y,t)\, t).
}
Thanks to Proposition 1 of \cite{sm_vr}, it holds that 
 $|x-y|^2 = \rho(x,y)^2 - \te(x,y)|x-y|^4$,
  where $\theta(x,y)>0$ is bounded on $M\x M$. 
  Applying Taylor's expansion,
  we  obtain that there exists $\eps_2>0$ and a bounded function $\td\Te: D_{\eps_2,\al} \to\Rnu$ so that
   on $D_{\eps_2,\al}$,  $\exp\{-\frac{\theta(x,y)|x-y|^4}{2t}\} = 1 + \td\Te(x,y,t)\, t^{4\al -1}$.
  Hence, by \rf{pq}, on $D_{\eps_2,\al}$ we have the relation
  \aaa{
  \lb{qe}
 q(x,y,t)  = \mc E(x,y,t)(1+\td\Te(x,y,t) \,t^{4\al-1}).
 }
  Next, By Theorem 3.22 of \cite{rosenberg}, for each  $k> \frac{d}2 + 2$,
\aaa{
\lb{hk}
h(x,y,t) = H_k(t,x,y) - Q_k \ast H_k(t,x,y),
}
where $H_k(t,x,y)$ is the heat kernel parametrix with the representation
\aaa{
\lb{parametrix}
H_k(t,x,y) = \eta(x,y) \mc E(x,y,t) (u_0(x,y) + \bar\Te(x,y,t) t),
}
and 
$Q_k \ast H_k(t,x,y) = \int_0^t d\te \int_M Q_k(\te,x, q) H_k(t-\te,q,y) \la_M(dq)$.
Function $Q_k: [0,T]\x M \x M\to \Rnu$ has an explicit representation for which
we refer the reader to \cite{rosenberg}, Section 3.2.1.
However, in our proof, we are only interested in estimate \rf{hk-rest} below 
which can be found in \cite{rosenberg} (proof of Proposition 3.23).
Namely, there exists a constant $C>0$ so that on $[0,T] \x M \x M$,
\aaa{
\lb{hk-rest}
|Q_k \ast H_k(t,x,y)|<C\,t^3.
}
Furthermore, in \rf{parametrix}, $\eta: M\x M \to [0,1]$ is a smooth function such that 
$\eta(x,y) = 1$ on $U_{\eps_3}(\diag(M\x M))$ for some $\eps_3>0$,
  and $\bar\Te(x,y,t)$ is bounded on $U_{\eps_3}(\diag(M\x M))\x [0,T]$
for some $T>0$. 
 Now \rf{hk}, \rf{parametrix}, and \rf{hk-rest} imply that there exists 
a bounded function $R:  M \x M \x [0,T]$ so that
  \aaa{
  \lb{hk1}
  h(x,y,t)=\eta(x,y)\mc E(x,y,t) (u_0(x,y) + \bar\Te(x,y,t) t) + R(x,y,t) t^3
  }
  for all $(x,y,t) \in M\x M \x [0,T]$.
  Moreover, the proof of Proposition 3.29 
  of \cite{rosenberg} implies that $u_0(x,x)=1$ and $\nab_M u_0(x,x)=0$. Therefore,
  applying Taylor's expansion to $u_0(x,y)$, we obtain that there exists $\eps_4>0$
  and a bounded function $\hat\Te: D_{\eps_4,\al}\to\Rnu$ so that on
  $D_{\eps_4,\al}$,
  \aaa{
  \lb{hk2}
   h(x,y,t)=\mc E(x,y,t) (1 + \hat\Te(x,y,t) t^{2\al}) + R(x,y,t) t^3.
   } 
   Relations \rf{qe} and \rf{hk2} imply that there exists $\eps>0$ so that
     \rf{2nd} holds on $D_{\eps,\al}$.
 Note that by construction, the bounds for $\Theta$ and $R$ do not depend on $\al$.
   \end{proof}
 \begin{lem}
   \lb{lem9}
Let $m,\beta>0$. The function $\Upsilon: (0,\infty) \to \Rnu, \quad  t\mto t^{-m} e^{-\frac1{2t^\beta}}$, 
is bounded by $\Big(2e^{-1} \frac{m}{\beta}\Big)^\frac{m}{\beta}$.
   \end{lem}
   \begin{proof}
   Computing the derivative in $t$, we obtain
   \aa{
 \frac{d}{dt}  \Big(t^{-m} e^{-\frac1{2t^\beta}} \Big) = e^{-\frac1{2t^\beta}} t^{-m-\beta -1}\Big(\frac{\beta}2 - m\, t^\beta\Big).
   }
   Taking into account that $\Upsilon$ tends to zero whenever $t\to 0$ or $t\to +\infty$, we conclude
   that the maximum of $\Upsilon$ is achieved at $t_0=\Big(\frac{\beta}{2m}\Big)^\frac1{\beta}$ and equals to
   $\Big(2e^{-1} \frac{m}{\beta}\Big)^\frac{m}{\beta}$. 
     \end{proof}
 \begin{lem}
 \label{lem2}
 Let $\al\in (\frac14, \frac12)$. Then, there exists $\eps>0$ and bounded
 functions $R_q$, $R_h$, $\Phi$, $\td \Phi$: $M\x M\x (0,\eps) \x (0,\eps) \to \Rnu$, so that
 for $x, z\in M$ and for $s,t \in (0,\eps)$,
\aaa{
& \int_M  q(x,y, t)  h(y,z,s)\la_M(dy)  =   h(x,z,t+s)(1 +  \Phi(x,z,t,s) t^{4\al-1}) \notag\\
 & + R_h(x,z,t, s) t^3,  \ \lb{78}\\
& \int_M  q(x,y, t)  q(y,z,s)\la_M(dy)  =  \int_M  q(x,y, t)  h(y,z,s)\la_M(dy) \notag \\
 &  \x (1 + \td \Phi(x,z,t,s) s^{4\al-1})+  R_q(x,z,t, s) \, s^3.  \lb{98}
 }
 Moreover, the functions $\Phi$ and $\td\Phi$  are bounded uniformly in $\al$, and
 the bound for
 $R_q$ and $R_h$ takes form \rf{bound1} with some constant $K>0${\rm\,:}
 \aaa{
 \lb{bound1}
 K \Big(\frac{d+3}{1-2\al}\Big)^\frac{d+3}{1-2\al}.
 }
 \end{lem}
    \begin{proof} 
 Let $U^x_{t}=\{y\in M: |y-x|<t^\al\}$, and let $V^x_{t} = M\dd U_{t}^x$.
 By Lemma \ref{lem9} (with $m=d+3$ and $\beta = 1-2\alpha$) and formula  \rf{pq}, there exists a constant $K$ so that
 \aaa{
 \lb{ineq}
\int_{V^x_{t}}q(x,y,t)h(y,z,s)\la_M(dy)  
<  K  t^{-d} e^{-\frac1{2t^{1-2\al}}} < K \Big(\frac{d+3}{1-2\al}\Big)^\frac{d+3}{1-2\al}  t^3.
}
In a similar manner, but with the use of formula \rf{hk1}, we obtain 
\aaa{
 \lb{ineq1}
 \int_{V^x_{t}}h(x,y,t)h(y,z,s)\la_M(dy) ) <   K \Big( \frac{d+3}{1-2\al}\Big)^\frac{d+3}{1-2\al}  t^3.
  }
  Without loss of generality, we  use the same constant $K$ as in \rf{ineq}. 
   Next, by Lemma \ref{lem1}, there exist functions 
$\Phi$, $\td R_h$: $M\x M \x (0,\eps) \x (0,\infty) \to \Rnu$,
so that
\mm{
 \int_{U_t^x} \! q(x,y,t)h(y,z,s)\la_M(dy) \\
   =  \int_{U_t^x} \! h(x,y,t)h(y,z,s)\la_M(dy) \big(1+ t^{4\al-1}  \Phi(x,z,t,s)\big)
+ \td R_h(x,z,t,s) \, t^3.
   }
   Moreover, by Lemma \ref{lem1}, $\Phi$ and $\td R_h$ are bounded uniformly in $\al \in (\frac14,\frac12)$.
By  \rf{ineq} and \rf{ineq1}, there exists a bounded
function $R_h: M\x M \x (0,\eps) \x (0,\infty) \to \Rnu$  so that \rf{78} holds.
The proof of \rf{98} follows the same steps.
 \end{proof}

\noindent 
Define the operator 
\aa{
Q_t: \C(M)\to \C(M), \, f \mto \int_M q(\,\cdot\,,y,t)f(y)\la_M(dy), \quad t>0.
}
 \begin{lem}
 \label{lem2prove}
 Let  the partition $\mc P$ of $[0,t]$ be such that
 $\tau = t_{n+1} - t_{n}$
 satisfies $\tau^{d+9} > |\mc P \dd \{t_{n+1}\}|$.
 Then, as $|\mc P| \to 0$, 
 \aaa{
 \label{2prove}
 (Q_{\Dl t_0} \cdots Q_{\Dl t_{n-1}} h\, (\, \cdot \, , y,
 \tau))(x)
 \to h(x,y,t)
 }
 uniformly in $(x,y,t) \in M\x M\x [a,b]$, $a>0$.
 \end{lem}
 \begin{proof}
 By the results of \cite{sm_vr} (p. 596, Corollary 2), we have
 \mm{
 \|(Q_{\Dl t_0} \cdots Q_{\Dl t_{n-1}} 
  -e^{-\frac{t-\tau}2 \lap_M}) \, p (\fdot, y, \tau)\| 
   \lt
   K \, t\, \|p ( \fdot  , y, \tau)\|_4 \: \sqrt{|\mc P \dd
  \{t_{n+1}\}|},
  }
 where  the norm $\|\fdot\|_4$ is defined in \cite{sm_vr} (see p. 593) as one of the equivalent norms in $\C^4(M)$. 
 Clearly,
 $\|p (\fdot, y, \tau)\|_4 < K_1 {\tau^{-(\frac{d}2+4)}}$, where $K_1>0$ is a constant.
 Since
 $|\mc P \dd \{t_{n+1}\}| < \tau^{d+9}$, we obtain that there exists a constant $K_2>0$ so that
 \aaa{
 \lb{es5}
 \|(Q_{\Dl t_0} \cdots Q_{\Dl t_{n-1}} - e^{-\frac{t-\tau}2 \lap_M})\, p (\, \cdot \, , y, \tau)\|
 < K_2 \, \sqrt{\tau} \to 0 \quad \text{as}\; |\mc P| \to 0.
 }
Next, by \rf{expansion},
 \mmm{
 \label{es7}
 (e^{-\frac{t-\tau}2 \lap_M}) \, p (\, \cdot \, , y, \tau )(x)= \int_M h(x,z,t-\tau)\, p(z,y,\tau)\la_M(dz) 
 = h(x,y, t-\tau) \\ -\frac{\tau}{2}\lap_M h(x,y,t-\tau)
 + \,\tau h(x,y, t-\tau)\Big(\frac16 \scal(x)
  +\frac1{16} \lap_M\lap_M\left.|x-\,\cdot\,|^2\right|_{x}\Big) \\ + \tau^{\frac32} R (x,t, \tau)
  \longrightarrow h(x,y,t), \quad \text{as} \; \tau \to 0,
  }
where the convergence holds uniformly in $(x,y,t) \in M \x M \x [a,b]$.
 This proves that as $|\mc P| \to 0$,
 $(Q_{\Dl t_0} \cdots Q_{\Dl t_{n-1}} p\, (\, \cdot \, , y, \tau))(x)
 \to h(x,y,t)$ uniformly in $(x,y,t) \in M \x M \x [a,b]$.
 Now \rf{2prove} follows from \rf{pq} and  \rf{98}.
 \end{proof}
 \begin{lem}
 \label{lem4}
 Let $\tau_i >0$ and $\sum_{i=1}^{n_\tau} \tau_i = \tau \to 0$ so that
 $\min_i \tau_i > \tau^\gm$ for some $\gm>1$.
 Then, there exists 
 $0<\eps<1$ such that $\sum_{i=1}^{n_\tau}  \tau_i^{1-\eps} \to 0$ as $\tau\to 0$.
  \end{lem}
 \begin{proof}
Pick $\eps<\gm^{-1}$. We obtain
$\sum_{i=1}^{n_\tau}  \tau_i^{1-\eps} \lt (\min_i \tau_i^\eps)^{-1} \tau < \tau^{1-\gm\eps} \to 0$.
\end{proof}
 \begin{proof}[Proof of Theorem~\ref{main}]
  Let $x,y \in M$. By Lemma \ref{lem2}, there exist  functions\\
  $\Te_1(x,y,\mc P)$,  $\Te_2(x,y,\mc P)$, $R_1(x,y,\mc P)$, and $R_2(x,y,\mc P)$ 
  bounded in $(x,y)\in M\x M$ and the points of the partition $\mc P$,  
  so that
 \mm{
 \hspace{-2mm} q_{\sss\mc P}(x, y, t)  
 = Q_{\Dl t_0} \cdots Q_{\Dl t_{n-2}}  h(\fdot,y, t_{n+1}-t_{n-1})(x)
  \big(1+ (\Dl t_{n})^{4\al-1}\Te_1(x,y, \mc P)\big) \\
 \x \big(1+ (\Dl t_{n-1})^{4\al-1}\Te_2(x,y, \mc P)\big) 
  + R_1(x,y,\mc P)(\Dl t_{n-1})^3 + R_2(x,y,\mc P)(\Dl t_{n})^3.
 }
 Let $N$ be the smallest number satisfying $\tau = t - t_{n-N} \gt |\mc
 P|^\frac1{d+9}$. In particular, this implies that $\tau - (t_{n-N+1}-t_{n-N}) <  |\mc
 P|^\frac1{d+9}$, and hence,
 $\tau < |\mc P|^\frac1{d+9} + |\mc P|$.
 Applying Lemma \ref{lem2} $(N+1)$ times, we obtain that there exist bounded functions $R_i(x,y,\mc P)$
 and $\Te_i(x,y,\mc P)$, 
 $i=0,\ldots,N$, so that
 \mmm{
 \lb{f90}
 q_{\sss \mc P}(x, y,t)  = \mc K_{t_{n-N}, \ldots, t_n} 
   \bigl(Q_{\Dl t_0} \cdots Q_{\Dl t_{n-N-1}}\, h(\, \cdot \, , y, t-t_{n-N})\bigr)(x) \\
   +  \mc R_{t_{n-N}, \ldots, t_n} ,
  } 
 where $\mc K_{t_{n-N}, \ldots, t_n}  = \prod_{i=0}^{N} \big(1+(\Dl t_{n-i})^{4\al-1}
 \Te_i(x,y,\mc P)\big)$, and
$\mc R_{t_{n-N}, \ldots, t_n} = \sum_{i=0}^{N} R_i(x,y,\mc P) (\Dl t_{n-i})^3$.
  By Lemma \ref{lem2prove} and the choice of $\tau$, as $|\mc P| \to 0$,
 \aaa{
  \label{lim11}
 \bigl(Q_{\Dl t_0} \cdots Q_{\Dl t_{n-N-1}}h\, (\, \cdot \, , y, t_{n}-t_{n-N})\bigr)(x) \to h(x,y,t)
 }
uniformly in $(x,y,t) \in M \x M \x [a,b]$.
It is clear that $\mc R_{t_{n-N}, \ldots, t_n}  \to 0$ uniformly in $(x,y,t) \in M \x M \x [a,b]$, since all $R_i$
are bounded
by the same constant which follows from the construction of $R_h$ in the proof of Lemma \ref{lem2}.
Let $\tau_i =\Dl  t_{n-i}$ for  $i=0,\ldots, N$.
To prove that $\mc K_{t_{n-N}, \ldots, t_n}  \to 1$, it suffices to show that
 $\sum_{i=0}^N\,\log(1+\tau_i^{4\al-1}\Te_i) \to 0$ as
 $|\mc P| \to 0$. Since
 $|\Te_i|$ are bounded by the same constant, as it is implied by the construction of $\Phi$
 in the proof of Lemma \ref{lem2},
we obtain that
\aaa{
\lb{in1}
 -2 \tau_i^{4\al-1}|\Te_i| < \log(1+\tau_i^{4\al-1}\Te_i) < \tau_i^{4\al-1}|\Te_i|
 }
when the mesh $|\mc P|$ is sufficiently small. Moreover, when $\tau$ and $|\mc P|$ are small,
$\tau^2 < \frac{\tau}2 < |\mc P|^{\frac1{d+9}}$. By the assumption of the theorem and the latter inequality,
\aaa{
\lb{in4}
 \min_i \tau_i > |\mc P|^L > \tau^{2L(d+9)}.
 }
Thus, we are in the conditions of Lemma \ref{lem4}. 
Pick $\al\in (\frac14, \frac12)$ sufficiently close to $\frac12$. 
By Lemma \ref{lem4},
 $\sum_{i=1}^N \tau_i^{4\al-1} \to 0$ as $|\mc P| \to 0$.
 The theorem is proved.
 \end{proof}
 \begin{cor}
Let the conditions of Theorem \ref{main} be satisfied, and
 let $u(x,t)$ be the
 solution of the Cauchy problem for \rf{heat} with the initial condition $u(0,x) = f(x)$, $f\in \C^2(M)$. Then 
 \aaa{
 \lb{cauchy}
 u(x,t) =\lim_{|\mc P|\to 0} \int_M q_{\sss \mc P}(x,y,t) f(y) \la_M(dy), 
 }
 where the limit is uniform in $(x,t)\in M\x [a,b]$, $a>0$. 
 \end{cor}
 \begin{proof}
 Since $u(x,t) = \int_M h(x,y,t) f(y) \la_M(dy)$, representation \rf{cauchy} immediately follows from Theorem \ref{main}.
 \end{proof}
\subsection{Rate of convergence}
Here we will estimate the rate of the convergence established in Theorem \ref{main}. Namely, we have the following result.
\begin{thm}[Uniform convergence rate]
\lb{rate}
Under the assumptions of Theorem \ref{main}, there exists a constant $K >0$ 
so that for all $(x,y,t) \in M\x M\x [a,b]$,
where $[a,b]\sub\Rnu$, $a>0$,
\aaa{
\lb{f10}
|h(x,y,t) - q_{\mc P}(x,y,t)|  < K |\mc P|^{\frac{1}{2(d+9)}}.
}
\end{thm}
\begin{proof}
In what follows, $\dl = 2-4\al$, 
$\Phi_i(x,y,\mc P)$, $i=1, \ldots, 8$, and $R_j(x,y,\mc P,\dl)$, $j=1,2$, will denote functions which are bounded in $(x,y)\in M\x M$ 
and $\mc P$ (for each fixed $\dl$), and 
the constant $K>0$ may differ from line to line, however we will use the same symbol for different constants.
 The bounds for $\Phi_i$ are uniform in $\dl$ (or in $\al$),
while the bounds for $R_j$ take form \rf{bound1}.
Also, we define $\gm = 2L(d+9)$, and without loss of generality assume that $\gm>1$.
As in the proof of Theorem \ref{main}, $\tau = t - t_{n-N}$.

From \rf{f90} it follows that 
\aaa{
\lb{e1}
\big| q_{\sss \mc P}(x, y,t)  - \mc K_{t_1, \ldots, t_{n-N}}
   \bigl(Q_{\Dl t_1} \cdots Q_{\Dl t_{n-N-1}}\, h(\, \cdot \, , y,\tau)\bigr)(x)\big| 
   \lt  K_\dl |\mc P|^2,
}
where $K_\dl = K \big(2(d+3)\dl^{-1}\big)^{2(d+3)\dl^{-1}}.$
Since  by the choice of $\tau$,
$|\mc P|^{\frac1{d+9}} \lt \tau < |\mc P| + |\mc P|^{\frac1{d+9}} < 2  |\mc P|^{\frac1{d+9}}$,
it suffices to find the convergence rate in the form $K\tau^\sg$ for some $\sg>0$.

First, we prove that there exists a function $\Phi$ bounded uniformly
in $x,y,\mc P$, and $\dl$ so that  
\aaa{
\lb{07.12.16}
\bigl(Q_{\Dl t_1} \cdots Q_{\Dl t_{n-N-1}}\, h(\, \cdot \, , y, \tau)\bigr)(x)  =
h(x,y,t) + \Phi(x,y,\mc P) \tau^\frac12.
}
By identity \rf{98} of Lemma \ref{lem2}, there exist functions 
$R_1(x,y,\mc P,\dl)$ and $\Phi_1(x,y,\mc P)$ bounded in $(x,y,\mc P)$ (for each fixed $\dl$) and such that
\mm{
 \bigl(Q_{\Dl t_1} \cdots Q_{\Dl t_{n-N-1}}\, q(\, \cdot \, , y, \tau)\bigr)(x)  \\  =
  (1+\Phi_1(x,y,\mc P)\tau^{1-\dl})\bigl(Q_{\Dl t_1} \cdots Q_{\Dl t_{n-N-1}}\, h(\, \cdot \, , y, \tau)\bigr)(x) + 
  R_1(x,y,\mc P, \dl)\tau^3.
}
Moreover, the bound for $\Phi_1$ is uniform in $\dl$, and the bound 
for $R_1$ is of form \rf{bound1}.
By the above relation and by \rf{pq}, one can find functions $\Phi_2(x,y,\mc P)$ and $R_2(x,y,\mc P,\dl)$ 
such that
\mmm{
\lb{fo9}
\bigl(Q_{\Dl t_1} \cdots Q_{\Dl t_{n-N-1}}\, h(\, \cdot \, , y, \tau)\bigr)(x) 
 = (1+\Phi_2(x,y,\mc P) \tau^{1-\dl})
 \\ \x \bigl(Q_{\Dl t_1} \cdots Q_{\Dl t_{n-N-1}}\, p(\, \cdot \, , y, \tau)\bigr)(x) 
+ R_2(x,y,\mc P, \dl)\tau^3. 
} 
Here the bound for $\Phi_2$ is uniform in $\dl$, while the bound for $R_2$ is of form \rf{bound1}.
By \rf{es5} and \rf{es7}, as well as by the boundedness of the heat kernel $h(x,y,t)$ and its derivatives
on $M\x M\x [a,b]$, $a>0$,
\mmm{
\lb{f8}
 \big(Q_{\Dl t_1} \cdots Q_{\Dl t_{n-N-1}} p (\, \cdot \, , y, \tau)\big)(x) =  \big(e^{-\frac{t-\tau}2 \lap_M} p (\, \cdot \, , y, \tau)\big)(x)
 + \Phi_3(x,y,\mc P) \tau^{\frac12}\\
= h(x,y,t-\tau) + \Phi_4(x,y,\mc P)\tau^\frac12 = h(x,y,t) +  \Phi_5(x,y,\mc P)\tau^\frac12.
 }
We remark that the bounds for $\Phi_i$, $i=1,\ldots, 5$,
do not depend on $\dl$ (or $\al$). Finally, substituting \rf{f8} into \rf{fo9} we obtain
\mm{
 \bigl(Q_{\Dl t_1} \cdots Q_{\Dl t_{n-N-1}}\, h(\, \cdot \, , y, \tau)\bigr)(x)   =
h(x,y,t) +  \Phi_6(x,y,\mc P)(\tau^\frac12 + \tau^{1-\dl}) \\ +  R_2(x,y,\mc P,\dl)\tau^3.
 }
Now let us estimate the term $\mc K_{t_1, \ldots, t_{n-N}}$. 
By \rf{in1}, \rf{in4}, and Lemma \ref{lem4},
\mmm{
\lb{in7}
-2K \tau^{1-\gm\dl} \lt -2  \sum\nolimits_{i=1}^N\, \tau_i^{4\al-1}|\Theta^{(i)}|  \lt \sum\nolimits_{i=1}^N\,
\log(1+\tau_i^{4\al-1}\Theta^{(i)}) \\ \lt \sum\nolimits_{i=1}^N \tau_i^{4\al-1}|\Theta^{(i)}|
\lt K \tau^{1-\gm\dl}.
}
Recall, that by Lemma \ref{lem2} and by the construction, the functions $\Theta^{(i)}(x,y,\mc P)$ are bounded in $x$, $y$, and $\mc P$ 
by a constant $K$ that does not depend on $i$ and $\dl$. 
Taking exponentials in \rf{in7}, we obtain
$e^{-2K \tau^{1-\gm\dl}} \lt \mc K_{t_1, \ldots, t_{n-N}}  \lt e^{K\tau^{1-\gm\dl}}$.
Note that for any $x>0$, $e^x-1 < x e^x$ and  $e^{-x} -1 > -x$. Therefore,
\aa{
- 2K \tau^{1-\gm\dl} <    \mc K_{t_1, \ldots, t_{n-N}} -1 < K \tau^{1-\gm\dl} e^{K\tau^{1-\gm\dl}}.
}
Hence, whenever $\tau<1$ and $\dl \lt \frac1{2\gm}$,  $\mc K_{t_1, \ldots, t_{n-N}}$ can be represented as
\aa{
\mc K_{t_1, \ldots, t_{n-N}}(x,y,\mc P) =1 +  \Phi_7(x,y,\mc P) \tau^{1-\gm\dl}.
}
Finally,
 \mmm{
\lb{07.12.16}
\mc K_{t_1, \ldots, t_{n-N}}(x,y,\mc P)\bigl(Q_{\Dl t_1} \cdots Q_{\Dl t_{n-N-1}}\, h(\, \cdot \, , y, \tau)\bigr)(x) =
h(x,y,t) \\ + \Phi_8(x,y,\mc P)(\tau^\frac12 + \tau^{1-\dl} + \tau^{1-\gm\dl}) + R_2(x,y,\mc P,\dl)\tau^3.
}
We remark, that the functions $\Phi_6,\Phi_7$, and $\Phi_8$ are bounded with the bounds not depending on $\dl$.
Note that the order of convergence in \rf{07.12.16} cannot be better than $\tau^\frac12$. Therefore, we can fix 
$\dl = \frac1{2\gm}$, or, which is the same, $\al = \frac12 - \frac1{8\gm}$ in both \rf{07.12.16} and \rf{bound1}. 
Hence,
by \rf{e1} and \rf{07.12.16}, 
\aaa{
\lb{final-id}
|q_{\sss \mc P}(x, y,t) - h(x,y,t)| < K \tau^{\frac12},
}
where $K$ is a constant.
Since $\tau < 2|\mc P|^{\frac1{d+9}}$, inequality \rf{final-id} implies \rf{f10}. 
\end{proof}
\begin{rem}
\rm As it is implied by the proofs of Lemmas \ref{lem1} and \ref{lem2}, the dependence of
the bounds (considered throughout all the proofs) on $\al$ comes only via formulas \rf{78} and \rf{98}.
We got rid of this dependence by the specific choice of $\al$.
The constant $K$  in Theorem \ref{rate} may also depend on
$a,b$, the bounds for $h$, $\frac{\pl}{\pl t} h$, $\lap_M h$ on $M\x M\x [a,b]$, as well as
the bounds of some smooth functions on $M$ such as $\scal(x)$, $\lap_M\lap_M|\fdot - y|^2_{x}$, and others.
\end{rem}
\subsection{Approximation of the heat kernel on an $n$-dimensional sphere}
We expect Theorem \ref{main} to be used by practitioners with the aid of a computer, since 
in most cases it is not possible
to obtain closed-form expressions for the integrals appearing in  \eqref{qp0}. Nevertheless, some procedures 
can be done to reduce the computational cost, and we will demonstrate one of them.

Let $S^d$ denote a unit sphere in $\Rnu^{d+1}$.
Our goal is to reduce the computational cost by ``removing" some integrals  in formula \eqref{qp0} that, in general, 
require to be evaluated numerically. 
We are going to compute the limit of
\aaa{
\lb{prod}
{\mb P}_n =  \prod_{i=0}^n \int_{S^d} p(x_{i},x_{i+1}, \Dl t_n) \la_{S^d}(dx_i)
}
in case of the uniform partition $\mc P$ of $[0,t]$. Here, $p(x_{i},x_{i+1}, \Dl t_n)$
is given by \rf{pdens},
$\Dl t_n = \frac{t}{n+1}$, $x_0 = x$, $x_{n+1} = y$.
This will allow us to simplify formula \rf{qp0} for $q_{\sss \mc P}$, where, roughly speaking,
the Brownian bridge density $q(x,y,t)$ will be substituted by
the Gaussian type density $p(x,y,t)$ given by \rf{pdens}, and also, to
provide a significant reduction on the computational cost by ``removing" $(n+1)$ numerical integrals from \eqref{qp0}. 
More precisely, we have the following result.
\begin{thm}
\lb{th3}
The heat kernel on $S^d$ has the following representation
\mm{
 h(x,y,t) =   e^{(\frac{d^2}{24} - \frac{d}6)t} 
  \lim\limits_{n\to\infty} (2\pi\Dl t_n)^{-\frac{(n+1)d}2} \\
 \x \int_{(S^d)^n} \exp\big\{ -\frac{\sum_{i=1}^{n+1} |x_i - x_{i-1}|^2}{2\Dl t_n}\big\} \la_{(S^d)^{n}}(dx_1 \ldots dx_n).
}
Moreover, the convergence is uniform in $(x,y,t)\in S^d \x S^d \x [a,b]$.
\end{thm}
\begin{proof}[Proof of Theorem \ref{th3}]
By applying asymptotic expansion \rf{expansion}, we obtain
\aaa{
\lb{expansion1}
\int_{S^d} \hspace{-1mm}  p(x,y, s) \la_{S^d}(dy)
   = 1 - s \big( \frac16 \scal(x)+  \frac1{16} \lap_{S^d}^2\left.|x-\,\cdot\,|^2\right|_{x}\big) \hspace{-1mm} +  s^{\frac32} \mc R(x,s).
}
Recall, that for a unit $d$-dimensional sphere, the scalar curvature $\scal(x)$ is constant and equals to $d(d-1)$.
Let us compute the bi-Laplacian $\lap_{S^d}^2 |x- y|^2$ with respect to $y$. Since $|x-y|^2 = 2 - 2(x,y)$, we start by computing
$\lap_{S^d} (x,y)$.  It is known that $\lap_{S^d} (x,y) = \left.-\lap\big[ (x, y|y|^{-1})\big]\right|_{|y|=1}$, where $\lap$ is the Laplacian in $\Rnu^{d+1}$.
It is immediate to compute $\lap |y|^{-1} = -(d-2)|y|^{-3}$ and $\nab |y|^{-1} = -y |y|^{-3}$, and therefore,
\aa{
\lap_{S^d} (x,y) = \hspace{-1.5mm}  \left. - \lap \big[|y|^{-1} (x,y)\big] \right|_{|y|=1} \hspace{-1mm} = \hspace{-1mm}\left. (x,y)\lap |y|^{-1} \hspace{-1mm} + 2 (x,\nab |y|^{-1})\right|_{|y|=1} \hspace{-1.4mm}  =  \hspace{-0.5mm} d \cdot \!(x,y).
}
Finally, 
\aa{
\left.\lap^2_{S^d}|x- y|^2\right|_{y=x} = \left. -2\, \lap^2_{S^d} (x,y)\right|_{y=x}
=  -2\, d^2 \cdot (x,y)|_{y=x} = -2  d^2.
}
Now \rf{expansion1} implies
\aa{
\int_{S^d} p(x,y,\Dl t_n) \la_{S^d}(dy)  = 1 + \Big(\frac{d}6  - \frac{d^2}{24}\Big) \Dl t_n +  (\Dl t_n)^{\frac32} \mc R(\Dl t_n),
}
where $\Dl t_n = \frac{t}{n+1}$. It is easy to see that the left-hand side of the above identity does not depend on $x$. Therefore,
the function $\mc R(x,\Dl t_n)$ in \rf{expansion1} does not depend on $x$ either, so in the above formula we write it as $\mc R(\Dl t_n)$.
Hence,
\aa{
\lim_{n\to \infty} P_n =  \lim_{n\to \infty} \big(1 + \Big(\frac{d}6  - \frac{d^2}{24}\Big) \Dl t_n +  (\Dl t_n)^{\frac32} \mc R(\Dl t_n)\big)^{n+1} = e^{(\frac{d}6  - \frac{d^2}{24})t}.
}
It remains to prove that this convergence is uniform in $(x,y,t)\in S^d\x S^d\x [a,b]$, $a>0$. 
Note that the product $\mb P_n$ does not depend on $x$ and $y$. Therefore, by Theorem \ref{main}, we just need to show that 
the convergence is uniform in $t\in [a,b]$, which follows from the fact that
the convergence of $(1+\frac{u}{n})^n$ to $e^u$ is uniform on compact sets. The theorem is proved.
\end{proof} 
\begin{rem}
\rm
The strategy used in the proof of Theorem \ref{th3} will work for any compact connected Riemannian manifold $M$ 
with the property that the value $\frac16 \scal(x)+  \frac1{16} \lap_{M}\lap_{M}\left.|x-\,\cdot\,|^2\right|_{x}$ is constant. 
\end{rem}

\end{document}